\documentclass[10pt,a4paper]{article}

\usepackage{graphicx}
\usepackage{amsmath}
\usepackage{amssymb}
\usepackage[numbers]{natbib}
\usepackage{microtype}
\usepackage[final]{showkeys}

\numberwithin{equation}{section}

\title{Fractional calculus modelling for unsteady unidirectional flow of incompressible fluids with
	time-dependent viscosity}

\author{$\text{Roberto Garra}_1$, $\text{Federico Polito}_2$\\
	\footnotesize (1) -- Dipartimento di Scienze di Base e Applicate per l'Ingegneria,\\ 
	\footnotesize Sapienza, Universit\`a di Roma,\\
	\footnotesize Via A. Scarpa 16, 00161 Rome, Italy.\\
	\footnotesize Email address: rolinipame@yahoo.it \\
	\footnotesize (2) -- Dipartimento di Matematica,
	\footnotesize Universit\`a di Torino,\\
	\footnotesize Via Carlo Alberto 10, 10123, Torino, Italy.\\
	\footnotesize Email address: federico.polito@unito.it
	}

\begin{document}

	\maketitle

	\begin{abstract}

		In this note we analyze a model for a unidirectional unsteady
		flow of a viscous incompressible fluid with time dependent viscosity.
		A possible way to take into account such behaviour is to introduce a memory formalism,
		including thus the time dependent viscosity by using an integro-differential term and therefore
		generalizing the classical equation of a Newtonian viscous fluid. A possible useful choice,
		in this framework, is to use a rheology based on stress/strain 
		relation generalized by fractional calculus modelling. This is a model that can be used in applied
		problems, taking into account a power law time variability of the viscosity coefficient.
		We find analytic solutions of initial value problems in an unbounded and bounded domain.
		Furthermore, we discuss the explicit solution in a meaningful particular case.

		\smallskip

	
	\end{abstract}

	\section{Introduction}  
		
		In recent years fractional calculus modelling has found many interesting
		applications in the description of viscoelasticity (see for example
		\citet{Mainardi}). It is well known that replacing the ordinary time derivative
		with the Caputo fractional derivative in the constitutive rheological equations
		we introduce a memory formalism in our problem. The Caputo
		fractional derivative is defined by a convolution of the ordinary integer-order derivative
		with a power law kernel. In this sense, it is evident the utility of a generalised
		rheology for complex media. Within this framework, several applications of fractional modelling
		of non-Newtonian fluids were recently investigated in many papers.
		For example \citet{Tong} treated the unsteady helical flow of a generalised
		Oldroyd-B fluid and \citet{Qi} studied the Stokes's first problem
		for such generalised Oldroyd-B fluid.
		Moreover, \citet{Vierua} studied the fractional Maxwell model
		and \citet{Shah} studied the fractional Burgers's fluid model.
		
		In this paper we discuss a model of unsteady unidirectional flow of an
		incompressible fluid with time-dependent viscosity in a simple geometry. 
		The main idea is that a possible way to take into account such variability
		is to generalise the equations of a viscous Newtonian fluid by introducing an
		integro-differential term, i.e.\ a memory formalism. Using this approach we include
		the time-dependent viscosity in the kernel of this 
		term thus suggesting an analytic way to understand the physical role played by the variable viscosity.
		
		A possible useful choice is to use a rheology based on a stress/strain 
		relation generalised by means of fractional calculus. We arrive then at
		a class of integro-differential equations based on a power law kernel.
		This analytic tool provides us useful mathematical instruments
		that can be used in applied problems which consider a viscosity coefficient
		with a power law time variability.
		
		We show that this fractional model of an unsteady flow
		is described by an integro-differential equation with a Riemann--Liouville
		fractional integral. We also discuss a general way to solve the initial value problem related
		to the master equation in an unbounded and bounded domain. 
		Finally, we discuss the explicit solution in an interesting case in order to better understand the
		physical meaning of our results.

	\section{Unsteady unidirectional flow of incompressible fluids with time-dependent
		viscosity by fractional calculus modelling}

		The classical stress/strain relations of a Newtonian viscous fluids are defined by:
	 	\begin{align}
	 		\label{stress}
			\tau_{ij} & = 2 \mu e_{ij}+\lambda \text{Tr}(e_{ij}) \\
			& =  2 \mu e_{ij}+\lambda (\nabla \cdot V), \notag
		\end{align}							
		where $\tau_{ij}$  is the Cauchy stress tensor, $\mu$ is dynamic viscosity of the medium,
		$e_{ij}$ is the strain tensor and $\lambda$ is the volume viscosity coefficient.
		As we are studying incompressible fluids, that is when $\nabla \cdot V =0$,
		then the second term in the right hand side of \eqref{stress} vanishes and we 
		have the simple relation
	  	\begin{equation}
	  		\label{strain}
			\tau_{ij}= 2 \mu e_{ij}.
		\end{equation}	
		In this note we suggest a way to take into account a possible
		time-dependent viscosity by generalising \eqref{strain} with an integro-differential
		relation:
		\begin{equation}
			\label{inte}
			\tau_{ij}= 2 \mu_0\int_0^t K(t-\tau)  e_{ij}(\tau) \, \mathrm d\tau, \qquad t >0,
		\end{equation}
		where $\mu_0$ is the constant viscosity of the specific memoryless case and 
		$K(\cdot)$ is a general kernel; its form depends on the type of time variability of the fluid viscosity.
		We notice that equation \eqref{inte} was first introduced by \citet{gerasimov} in relation to the application
		of generalised laws of deformation to problems of internal friction (see also \citet{Mainardi} and
		\citet{uchaikin}, pages 285--286).
		
		This proposal has a clear physical meaning: the role played by the viscosity time
		dependence is to introduce a memory formalism and hereditary effects on the fluid flow.
						
		We now propose a useful way to generalise \eqref{strain}, specialising \eqref{inte} with a power law kernel.
		This generalisation is based on the application of the Riemann--Liouville fractional integral. Hence, we have
	 	\begin{equation}
	 		\label{gen}								
			\tau_{ij}= 2 \mu_0 J_t^{\alpha} e_{ij}, \qquad \alpha \in (0,1), \: t>0,
		\end{equation}
		where $J_t^{\alpha}$ is the Riemann--Liouville time fractional integral defined by \citep{Podl}
		\begin{equation}
			J^{\beta}_t f(t) = \frac{1}{\Gamma(\beta)}\int_0^{t}(t-\tau)^{\beta-1}f(\tau)
			\, \mathrm d\tau, \qquad \beta > 0, \: t>0.
		\end{equation}
		Within our framework, the meaning of \eqref{gen} is the following: through fractional
		integration we take into account a power law time variability 
		of the viscosity coefficient i.e.\ we include non-local effects. This can be interesting
		in relation to the study in the field of time variable viscosity in general (for example
		in solid suspensions \citep{alloy}).
		Here we recall that there are many different approaches to fractional rheologies starting from the seminal works
		of Scott--Blair and Gerasimov (for a historical review the reader can refer to \citet{main}).
		However we are not discussing the rheology of a specific class of viscoelastic media, we
		used the fractional formalism only in order to include the time variability of the viscosity
		by means of convolution.
		Therefore our model cannot be directly related to the already classical fractional rheologies
		widely discussed in the literature (see for a recent review \citet{surguladze}).
		Instead we are going to study the flow of a viscous fluid with a power-law time dependent
		viscosity. In general it leads to an integro-differential generalisation of the Navier--Stokes
		equation. In the following we will discuss a particular case as the more general case leads to
		a complex nonlinear integro-differential equation which however could be the object of
		a further investigation.
		
		The equations of flow of an incompressible fluid without volume force
		in Cartesian coordinates are
		\begin{align}
			\begin{cases}
				\rho\frac{\mathrm D V}{\mathrm Dt}= -\nabla p+ \nabla\cdot \tau, \\
				\nabla \cdot V = 0, 
			\end{cases}
		\end{align}							
		where $\mathrm D / \mathrm Dt = \partial_t+ (V \cdot \nabla)$ is the material derivative.
		For simplicity, we assume a constant density $\rho = 1$.
		We treat the unidirectional motion of this fluid on a plate.
		The flow of a unidirectional velocity field, in Cartesian coordinates is given by
		\begin{equation}
			V = u(y,t)\hat{\imath},
		\end{equation}
		where $u$ is the velocity in the $x$-direction, $\hat{\imath}$ is the unit vector in
		the $x$-direction, $x$ is the coordinate along the plate
		and $y$ is the coordinate normal to the plate.
		Considering these hypotheses on the velocity field, taking $\nabla p = 0$ (constant pressure),
		and using \eqref{gen} we have that
 		\begin{equation}
 			\label{u}
			\partial_t u = \mu_0 \partial_y(J_t^{\alpha}\partial_y u), \qquad t \ge0.
		\end{equation}
		This is an integro-differential equation with a fractional Riemann--Liouville integral.
		
		In the next, we first
		discuss the analytic solution of the following initial value problem (IVP) in an unbounded domain.
		The following step will be, in Section \ref{sect} the more realistic case of a bounded domain.
		\begin{align}
			\begin{cases}
				\partial_t u = \mu_0 \partial_y(J_t^{\alpha}\partial_y u), & y \in \mathbb{R},\;t \ge 0,\\
				u(y,0)= f(y).
			\end{cases}
		\end{align}
		In order to solve this IVP we make use of the Laplace--Fourier transform.
		Recalling that the Laplace tranform of the fractional integral is given by \citep{Podl}
		\begin{equation}
			\mathcal{L}(J_t^{\alpha}u(y,t))(s)= s^{-\alpha}\tilde{u}(y,s),
		\end{equation}
		where $s$ is the Laplace transform parameter and $\tilde{u}(y,s)$ is the Laplace
		transform of $u(y,t)$, we immediately obtain
		\begin{equation} 
			s\hat{\tilde{u}}(\omega,s)-\hat{f}(\omega)= -\omega^2\mu_0 s^{-\alpha}\hat{\tilde{u}}(\omega,s),
		\end{equation}
		where $\hat{\tilde{u}}(\omega,s)$ is the Laplace--Fourier tranform of $u(y,t)$. Therefore
		\begin{equation}
			\label{l-f}
			\hat{\tilde{u}}(\omega,s)=\frac{s^{\alpha}\hat{f}(\omega)}{s^{\alpha+1}+\mu_0\omega^2}.
		\end{equation}
		We now recall the following useful relation \citep{Podl}:
		\begin{equation}
			\int_0^\infty e^{-st} t^{\mu -1} E_{\nu,\mu} (a t^\nu) \mathrm dt = \frac{s^{\nu-\mu}}{s^\nu - a},
		\end{equation}
		where
		\begin{equation}
			E_{\nu,\mu}(z) = \sum_{r=0}^\infty \frac{z^r}{\Gamma(\nu r+ \mu)}, \qquad z \in \mathbb{R},
		\end{equation}
		is the Mittag--Leffler function. The Laplace transform in \eqref{l-f} can be easily
		inverted now, leading to the Fourier transform $\hat{u}(\omega,t)$:
		\begin{equation}
			\hat{u}(\omega,t) = \hat{f}(\omega) E_{\alpha+1,1}(-\mu_0 \omega^2 t^{\alpha+1}).
		\end{equation}
		By recurring to relations (14.7) and (14.8) of \citet{hau} we can invert also the Fourier
		transform. With this inversion we can write the solution $u(y,t)$ by means of a
		Fox function. Indeed we have
		\begin{align}
			u(y,t) = {} & f(y) \, \ast \, \left( \frac{1}{|y|} H^{1,0}_{1,1} \left[ \left.
			\frac{|y|^2}{\mu_0 t^{\alpha+1}} \right|
			\begin{array}{l}
				(1,\alpha+1) \\
				(1,2)
			\end{array}
			\right] \right),
			\qquad t \ge 0, \: y \in \mathbb{R},
		\end{align}
		where $\ast$ is the Fourier convolution.
		By means of simple manipulations and considering the representation
		of the Fox function as a contour integral \citep{hau} we obtain that
		\begin{align}
			\label{acqua}
			u(y,t) = {} & f(y) \, \ast \, \left( \frac{1}{|y|} \frac{1}{2 \pi i} \int_C \frac{\Gamma(1+2z)}{
			\Gamma(1+(\alpha+1)z)} \left( \frac{|y|^2}{\mu_0 t^{\alpha+1}}
			\right)^{-z} \mathrm d z \right) \\
			= {} & f(y) \, \ast \, \left( \frac{1}{\mu_0^{1/2} t^{\frac{\alpha+1}{2}}} \frac{1}{2\pi i}
			\int_C \frac{\Gamma(1+2z)}{\Gamma(1+(\alpha+1)z)}
			\left( \frac{|y|}{\mu_0^{1/2} t^{\frac{\alpha+1}{2}}} \right)^{-2z-1} \mathrm d z \right) \notag \\
			= {} & f(y) \, \ast \, \left( \frac{1}{2\mu_0^{1/2} t^{\frac{\alpha+1}{2}}} \frac{1}{2\pi i}
			\int_C \frac{\Gamma(\tau)}{\Gamma\left( 1+(\alpha+1)\frac{\tau-1}{2} \right)}
			\left( \frac{|y|}{\mu_0^{1/2} t^{\frac{\alpha+1}{2}}} \right)^{-\tau} \mathrm d \tau \right) \notag \\
			= {} & f(y) \, \ast \, \left( \frac{1}{2\mu_0^{1/2} t^{\frac{\alpha+1}{2}}} \frac{1}{2\pi i}
			\int_C \frac{\Gamma(\tau)}{\Gamma\left( 1-\frac{\alpha+1}{2} + \frac{\alpha+1}{2} \tau \right)}
			\left( \frac{|y|}{\mu_0^{1/2} t^{\frac{\alpha+1}{2}}} \right)^{-\tau} \mathrm d \tau \right). \notag 
		\end{align}
		In the last step of \eqref{acqua} we recognize the Mellin--Barnes representation
		of the Wright special function, defined also as		
		\begin{equation}
			W_{\nu,\mu}(z) = \sum_{r=0}^\infty \frac{z^r}{r!\Gamma(\nu r+\mu)},
			\quad \nu>-1, \: \mu \in \mathbb{R}, \: z \in \mathbb{C}.
		\end{equation}
		We therefore arrive at the final following form.
		\begin{align}
			\label{treno}
			u(y,t) = {} & f(y) \, \ast \, \left( \frac{1}{2\mu_0^{1/2} t^{\frac{\alpha+1}{2}}} W_{-\frac{\alpha+1}{2},
			1-\frac{\alpha+1}{2}} \left( -\frac{|y|}{\mu_0^{1/2} t^{\frac{\alpha+1}{2}}} \right) \right),
			\quad t \ge 0, \: y \in \mathbb{R}.
		\end{align}
		
		Note that the above solution coincides with the solution to the superdiffusive fractional
		diffusion equation with the initial condition $\partial_t u(y,t)|_{t=0}=0$.
		Note also that the well-known fundamental solution to the fractional diffusion equation
		written in terms of Wright functions
		(see e.g.\ \citet{Mainardi1}) can be retrieved by letting $f(y) = \delta(y)$.		
		
		This result is apparently not surprising, but it forces us to make some considerations
		on a sensitive topic related to fractional calculus, i.e.\
		the validity of the semigroup property for fractional derivatives.
		Actually, going back to the master equation \eqref{u} and applying the Caputo fractional
		derivative to both sides, we obtain
		\begin{equation}
			\label{gr}
			\partial^{\alpha}_t\partial_t u = \mu_0 \partial_{yy} u, \qquad t\ge 0, \: y \in \mathbb{R}.
		\end{equation} 
		We recall that the definition of the Caputo partial fractional derivative $\partial^{\alpha}_t$,
		$\alpha \in (0,1)$, is
		\citep{Podl}
		\begin{align}
			\partial_t^{\alpha}g(y,t) & = J_t^{1-\alpha}\partial_t g(y,t) \\
			& = \frac{1}{\Gamma(1-\alpha)}\int_0^{t}(t-\tau)^{\alpha}\partial_\tau g(y,\tau)\, \mathrm d\tau. \notag
		\end{align}
		Equation \eqref{gr} is a sequential fractional differential equation,
		different from the fractional superdiffusive equation
		thoroughly treated in literature.
		The reason of the difference is that the semigroup property is not always satisfied by Caputo fractional 
		derivatives, i.e. $\partial_t^{\alpha}\partial_t^{\beta} u \neq \partial_t^{\alpha+\beta} u
		\ne \partial_t^{\beta}\partial_t^{\alpha} u$.
		
		In literature, several attempts to solve equations with fractional
		sequential derivatives as in \eqref{gr} are present (see for example \citet{Podl}),
		but they require initial conditions with a not clear physical meaning.
		However, in our case we have a well posed problem because
		the master equation is basically a first order integro-differential equation.
		The difference with the superdiffusive fractional case is
		in the initial conditions. In order to solve \eqref{u} it is sufficient one initial condition,
		while in the superdiffusive case two conditions are mandatory.
		We remark, however, that the analytic solution of \eqref{u} coincides with that of the superdiffusive
		fractional equation in the specific case $\partial_t u(y,t)|_{t=0}=0$.
		 	
		In the next section we outline a general way to solve analytically the master equation in a bounded domain.
		Moreover, we discuss the exact solution for an interesting application which allows us
		to highlight the physical meaning for this model. 

	\section{Exact solution to the master equation in a bounded domain}
	
		\label{sect}
		A number of papers are devoted to the analytic investigation of exact solutions to
		fractional integro-differential equations in bounded domains.
		Here we recall for example that \citet{Agrawal} found  a general way to solve
		the fractional diffusion equation
		in a bounded domain by using Laplace and finite sine tranform techniques.
	
		Here we are going to make use of similar methods to solve the following initial/boundary value problem:
		\begin{align}
			\label{mast1}
			\begin{cases}
				\partial_t u= \mu_0 J^{\alpha}_t \partial_{yy}u, & \alpha \in (0,1), \\
				u(0,t)=u(h,t) = 0, & t\geq 0, \\
				u(y,0)=f(y), & 0 \le y \le h .
			\end{cases}
		\end{align}
		The finite sine transform of \eqref{mast1} is simply given by
		\begin{equation}
			\partial_t\bar{u}(n,t) = \mu_0 a^2 n^2 J^{\alpha}_t\bar{u}(n,t),
		\end{equation}
		where $n$ is a wave-number and $a = \pi/h$. Note also that $\bar{u}(n,t)$ is the finite sine transform of $u(y,t)$
		and is defined as
		\begin{equation}
			\bar{u}(n,t)=\int_0^h u(y,t)\sin(a\,n\,y)\, \mathrm dy.
		\end{equation}
		Therefore, the related initial value problem (in wave-number domain) is
		\begin{align}
			\label{sin}
			\begin{cases}
				\partial_t\bar{u}(n,t) = \mu_0 a^2 n^2 J^{\alpha}_t\bar{u}(n,t),\\
				\bar{u}(n,0)=\int_0^h f(y) \sin(a\,n\,y)\, \mathrm dy.
			\end{cases}
		\end{align}
		By taking the Laplace transform of \eqref{sin}, we obtain
		\begin{equation}
			\tilde{\bar{u}}=\frac{s^{\alpha}\bar{u}(n,0)}{s^{\alpha+1}+\mu_0(an)^2}.
		\end{equation}
		Finally we can write that
		\begin{equation}
			\label{sol}
			u(y,t)= \frac{2}{h}\sum_{n=1}^{\infty}
			E_{\alpha+1}\left(-\sqrt{\mu_0}\,a^2n^2t^{\alpha+1}\right)
			\sin(a\,n\,y)\int_{0}^{h}f(z)\sin(a\,n\,z) \mathrm dz,
		\end{equation}
		where $t \ge 0$, $y \in [0,h]$.
		
		This is a general way to solve the fractional diffusion equation in a bounded domain, with the
		additional constraint of
		homogeneous boundary conditions.
		
\begin{figure}
			\centering
			\includegraphics[scale=.8]{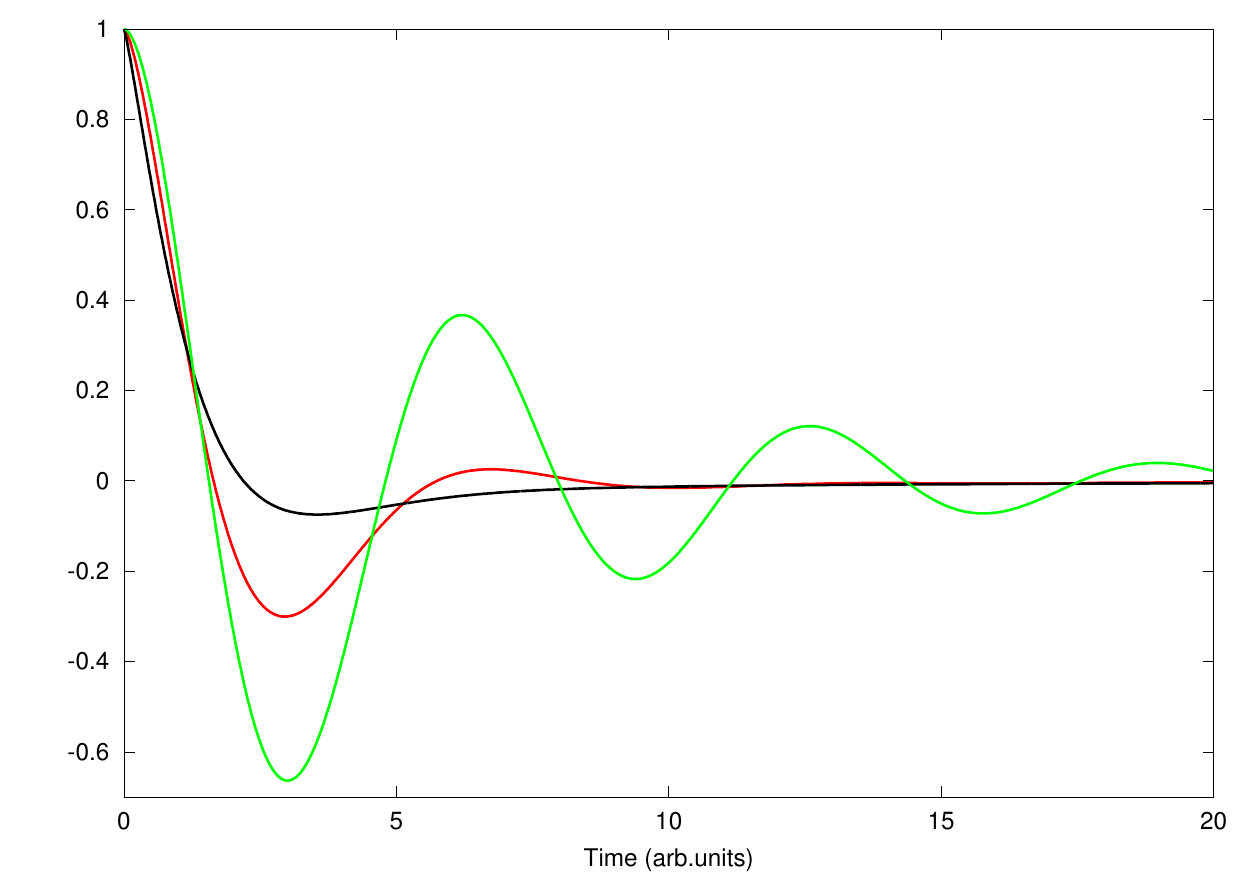}
			\caption{\label{fig}The Mittag--Leffler function $E_{\gamma,1}(-t^\gamma)$ for
				$\gamma=1.2$ (black), $\gamma=1.5$
				(red) and $\gamma=1.8$ (green).}
		\end{figure}

		Now we discuss an application of this analytic way to solve \eqref{mast1} in a simple case.
		We study the following problem:
		\begin{align}
			\begin{cases}
				\partial_t u= \mu_0 \partial_{yy}J^{\alpha}_t u, & \alpha \in (0,1),\\
				u(0,t)=u(h,t) = 0, & t\geq 0,\\
				u(y,0)=U_0= \text{const}, & 0 \le y \le h. 
			\end{cases}
		\end{align}
		The physical meaning of these conditions is clear: we have an initial constant velocity field
		$u(y,0)=U_0$ in a confined plane, with
		null boundary conditions at both sides.
		By using \eqref{sol}, we find
		\begin{equation}
			u(y,t)= U_0\frac{2}{\pi}\sum_{n= 1}^{\infty}
			E_{\alpha+1}\left(-\sqrt{\mu_0}\, a^2n^2t^{\alpha+1}\right)\frac{\sin(a\,n\,y)}{n}
			\left(1-\cos(n\pi)\right).
		\end{equation} 
		We have to notice that the Mittag--Leffler function
		of order $1<(\alpha+1)<2$ is also called in literature as a 
		fractional oscillation (see Fig.\ \ref{fig}).
		The result of memory effects due to time variability of
		the viscosity coefficient is thus an anomalous fractional oscillation. Furthermore,
		we can recover the exact periodic solution for $\alpha=1$. 
		As expected, this is the solution that can be found by separation of variables from the ordinary wave equation.
		
		Moreover, this solution suggests that fractional modelling is a good instrument 
		to describe global damping due to local time variability of the viscosity coefficient.
 	
	\section{Conclusions}
	
		In this note we suggested the use of a generalised stress/strain relation to consider a
		time variability of the fluid viscosity coefficient.
		Treating a unidimensional model, we showed that, under reasonable mathematical hypotheses,
		the fluid flow is described by a linear integro-differential equation
		involving a Riemann--Liouville fractional integral. We also discussed, by means of
		an exactly solvable example,
		the meaning and utility of this model also in relation to the study
		of the behaviour of complex fluids with time variable viscosity.
		Moreover, we notice that the real order of integration is a free parameter that can be
		exploited in order to have a better agreement with the experimental data.
		It is also relevant to understand the physical meaning of this real
		order of integration in relation to the characteristics of the fluid.
		Finally, we remark that within this framework we take into account
		possible emerging empirical power law behaviours by means of simple analytic models. 

	\bibliography{non3} 
	\bibliographystyle{unsrtnat}
	\nocite{*}

\end{document}